\documentclass[10pt]{article}
\usepackage{amsmath, amsfonts, amsthm, amssymb,verbatim}
        \newtheorem{theorem}{Theorem}[section]

\numberwithin{equation}{section}

\def\og{\leavevmode\raise.3ex\hbox{$\scriptscriptstyle\langle\!\langle$~}}
\def\fg{\leavevmode\raise.3ex\hbox{~$\!\scriptscriptstyle\,\rangle\!\rangle$}}

\newcommand \lam    {\lambda} 
\newcommand \RR     {{\, R\!\!\!\!\!\!I~ \,}}
\newcommand \del    {{\partial}}
\newcommand \eps    {\varepsilon}
\newcommand \RN     {\RR^N} 
     
\newcommand \CC     {\mathcal C}
  
\newcommand \WW     {\mathcal W} 
\newcommand \hatr   {\widehat r} 

\newcommand \be        {\begin{equation}}
\newcommand \ee        {\end{equation}}
\newcommand \Bzero  {{\mathcal B}_{\delta_0}}
\newcommand \Bone  {{\mathcal B}_{\delta_1}}
\newcommand \bLam {\underline \Lambda}
\newcommand \Lamb {\overline \Lambda}

\newcommand \hatl   {\widehat l}

\newcommand \bm     {{\underline m}} 
\newcommand \mb     {{\overline m}}  

\newcommand \Id 	{\mbox{Id}}

\begin{document}

\title{Singular limits for the Riemann problem. General diffusion, relaxation, and boundary conditions}
\author{K.T. Joseph$^1$ and Philippe G. LeFloch$^2$}

\date{February, 2006}

\maketitle

\footnotetext[1]{School of Mathematics, Tata Institute of Fundamental Research, Homi Bhabha Road,
Mumbai 400005, India. E-mail: ktj@math.tifr.res.in}

\footnotetext[2]{Laboratoire Jacques-Louis Lions \& Centre National de la Recherche Scientifique (CNRS),
Universit\'e Pierre et Marie Curie (Paris 6), 4 Place Jussieu,  75252 Paris, France.
E-mail: pgLeFloch@gmail.com. {\tt To cite this article:} C.R. Math. Acad. Sci. Paris 344 (2007), 59--64.   
}

\begin{abstract} 
We consider self-similar approximations of non-linear hyperbolic systems in one space
dimension with Riemann initial data, especially the system $\del_t u_\eps + A(u_\eps) \, \del_x
u_\eps = \eps \, t \, \del_x ( B(u_\eps) \, \del_x u_\eps)$, with $\eps >0$. We assume that the matrix 
$A(u)$ is strictly hyperbolic and that the diffusion matrix satisfies $|B(u) - \Id|<<1$. 
No genuine non-linearity assumption is required. We show the existence of a smooth, self-similar solution
$u_\eps = u_\eps(x/t)$ which has bounded total variation, uniformly in the diffusion parameter $\eps>0$. 
In the limit $\eps \rightarrow 0$, the functions $u_\eps$ converge towards a solution of the Riemann problem
associated with the hyperbolic system. A similar result is established for the relaxation approximation
$\del_t u^\eps + \del_x v^\eps = 0$, 
$\del_t v^\eps + a^2 \, B(u) \del_x u^\eps  = ( f(u^\eps) - v^\eps )/(\eps \, t)$. 
We also cover the boundary-value problem in a half-space for the same regularizations.
\end{abstract} 
 

\section{Introduction}
\label{Introd}

In this note we present a continuation of the work by the authors 
\cite{JL2}--\cite{JL5} on self-similar regularizations of the  
Riemann problem associated with a non-linear strictly hyperbolic system in one-space dimension:
\be
\del_t u + A(u) \, \del_x u = 0, \qquad u=u(t,x) \in \Bzero, \quad 
t>0, \, x \in \RR,  
\label{1.1}
\ee
with piecewise constant, initial data 
\be
u(0,x) = u_l  \,  \mbox{ for } x < 0; \quad u_r  \, \mbox{ for } x >0, 
\label{1.2}
\ee
where $u_l, u_r$ are constant states in $\Bzero$.  Here, $\Bzero \subset \RN$ denotes the ball centered at the origin and 
with radius $\delta_0>0$, and, for all $u \in \Bzero$,
$A(u)$ is assumed to have distinct and real eigenvalues
$\lam_1(u) < \ldots < \lam_N(u)$ and basis of left- and right-eigenvectors $l_j(u), r_j(u)$, ($1 \leq j \leq N$). 
 
Following Dafermos \cite{Dafermos1,Dafermos2} and Slemrod \cite{Slemrod} who advocate the use 
of self-similar regularizations to capture the whole wave fan structure of the Riemann problem, 
we consider solutions constructed by self-similar vanishing diffusion 
associated with a {\sl general} diffusion matrix $B= B(u)$, that is we search for solutions of 
\be
\del_t u_\eps + A(u_\eps) \, \del_x u_\eps 
= 
\eps \, t \, \del_x \bigl( B(u_\eps) \, \del_x u_\eps\bigr),
\quad \eps >0.     
\label{1.3}
\ee 
Due to the choice of the scaling $\eps \, t$, this system admits solutions $u_\eps = u_\eps(x/t)$, 
and, therefore, we refer to (\ref{1.2})-(\ref{1.3}) as the {\sl self-similar diffusive Riemann problem.} 
The matrix $B=B(u)$ is assumed to depend smoothly upon $u$ and to remain sufficient close to the identity matrix,
that is,  for some given matrix norm and for $\eta>0$ sufficiently small 
\be
\sup_{u \in \Bzero} |B(u) - \Id | \leq \eta. 
\label{1.4}
\ee
The method of analysis introduced below is not a~priori restricted to (\ref{1.3})-(\ref{1.4}), 
and generalizations are discussed at the end of this note and in \cite{JL5,JL6}. 

The techniques developed so far for general hyperbolic systems 
(see \cite{Tzavaras,LeFlochTzavaras} and \cite{JL2}--\cite{JL4}) 
were restricted to regularizations based on the identity diffusion matrix. The new approach introduced here 
allows us to cover classes of approximations based on general diffusion matrices (or relaxation terms, see below). 
This degree of generality is especially important for non-conservative systems  \cite{LeFloch1,LeFloch2}  
and for the boundary-value problem  \cite{JL2}, whose solutions are known to {\sl strongly depend} 
upon the specific regularization.

\section{Main results} 

By the property of propagation at finite speed, a self-similar solution $u=u(y)$ of the Riemann problem
is constant outside a sufficiently large, compact interval $[-L,L]$, i.e.: $u(y) = u_l$ for $y <-L$
and $u(y) = u_r$ for $y >L$. 
As is customary, we assume that $\delta_0$ is sufficiently small so that   
the wave speeds $\lam_j(u)$ remain close to the constant speeds $\lam_j(0)$ and
are uniformly separated in the sense that
$$
\bLam_j \leq \lam_j(u) \leq \Lamb_j, \quad u \in \Bzero,   
$$  
for some constants $-L < \bLam_1 < \Lamb_1 < \bLam_2 < \ldots < \bLam_N < \Lamb_N <L$. 

\

We establish the following theorem. 

\

\begin{theorem} 
\label{1-1} 
Consider the non-linear, strictly hyperbolic system (\ref{1.1}) together with its parabolic regularization (\ref{1.3})-(\ref{1.4}). 
There exist (sufficiently small) constants $\delta_1, \eta>0$ 
and a (sufficiently large) constant $C_0 > 0$ such that 
for any initial data $u_l, u_r \in \Bone$ the self-similar diffusive Riemann problem (\ref{1.2})-(\ref{1.3})  
admits a smooth solution $u^\eps=u^\eps(x/t) \in \Bzero$ defined for all $y=x/t \in [-L,L]$,
which has uniformly bounded variation,
$$
TV_{-L}^L(u_\eps) \leq C_0 \, |u_r - u_l|, 
$$ 
and converges strongly to some limit $u:[-L,L] \to \Bzero$: 
$$
u^\eps \to u \, \mbox{ in the $L^1$ norm, as } \eps \to 0. 
$$  
The limit function satisfies the following properties. The function $y \mapsto u(y)$ has bounded total variation, 
that is, $TV_{-L}^L(u) \leq C_0 \, |u_r - u_l|$, 
and is constant on each interval $[\Lamb_j, \bLam_{j+1}]$. If (\ref{1.1}) is a system of conservation laws,  
i.e.~$A=Df$ for some flux $f:\Bzero \to \RN$, then the limit is a distributional solution of 
\be
\del_t u + \del_x f(u) = 0.
\label{1.5}
\ee 
If $U,F): \Bone \to \RR \times\RN$ is an entropy / entropy flux pair associated with (\ref{1.5}) 
and the diffusion matrix satisfies the convexity-like condition $\nabla^2 
U \cdot B \geq 0$,
then the solution $u$ satisfies the entropy inequality 
\be
\del_t U(u) + \del_x F(u) \leq 0. 
\label{1.6}
\ee 
\end{theorem} 

\

We have also the following description of the wave curves. 

\

\begin{theorem}
\label{1-2} 
With the notation and assumptions in Theorem~\ref{1-1}, 
to each $j$-characteristic family and each left-hand state $u_l$ one can associate 
a {\rm $j$-wave curve} 
$$ 
\WW_j(u_l) := \bigl\{ u_r = \psi_j(m; u_l) \, / \, m \in (\bm_j, \mb_j) \bigr\}, 
$$   
issuing from $u_l$, which, by definition, is made of all right-hand states $u_r$ attainable by a Riemann solution $u=u(y)$, 
with left-hand state $u_l$, by using only $j$-waves, that is such that 
$$
u(y) = u_l \mbox{ for } y <  \bLam_j; \qquad u_r \mbox{ for } y > \Lamb_j. 
$$ 
Moreover, the mapping $\psi_j : (\bm_j, \mb_j) \times \Bone \to \Bzero$ 
is Lipschitz continuous with respect to both arguments, and for some small constant $c>0$  
$$
\del_m \psi_j(m; u_l) \in \CC_j 
:= \bigl\{ w \in \RN \, / \, \bigl| w \cdot l_j(0) \bigr| \geq (1-c) \, |w | \bigr\}. 
$$ 
Moreover, the characteristic component
$y \mapsto \alpha_j(y) := l_j(0) \cdot u'(y)$ 
is a non-negative measure in the interval $[\bLam_j, \Lamb_j]$.  
\end{theorem} 

\

For the proof of these results as well as a characterization of the limit when (\ref{1.1}) is a general non-conservative 
system, we refer to \cite{JL5,JL6}. 

We will here sketch the proof of Theorem~\ref{1-1}. To handle general diffusion matrix $B(u)$, the following 
{\sl generalized eigenvalue problem} is introduced 
$$
\bigl( - y + A(u) \bigr) \, \hatr_j(u,y) = \mu_j(u,y) \, B(u) \,
\hatr_j(u,y),
$$
$$
\hatl_j(u,y) \cdot \bigl( - y + A(u) \bigr) = \mu_j(u,y) \, \hatl_j(u,y) \cdot B(u). 
$$ 
In view of (\ref{1.4}), one has $\hatr_j(u,y) = r_j(u) + O(\eta)$ and $\hatl_j(u,y) = l_j(u) + O(\eta)$. 
The proof relies on a suitable {\sl asymptotic expansion}  
of the solution $u_\eps = u_\eps(x/t)$, of the form 
$$
u_\eps' = \sum_j  \, a_j^\eps \, \hatr_j(u_\eps,\cdot) 
\quad \mbox{ with } \, a_j^\eps  := \hatl_j(u_\eps, \cdot) \cdot u_\eps'. 
$$ 
Omitting $\eps$, we deduce that 
the components $a_j$ satisfy a {\sl coupled system of $N$ differential equations:} 
$$ 
a_i' - {\mu_i(u, \cdot) \over \eps} \, a_i 
+ \sum_j \pi_{ij}(u, \cdot)\, a_j 
= 
Q_i(u; \cdot) := \sum_{j, k} \kappa_{ijk}(u, \cdot) \, a_j \, a_k,    
$$
where 
$$ 
\pi_{ij}(u, \cdot) 
   : =   \hatl_i(u, \cdot) \cdot B(u) \,  \del_y \hatr_j(u,\cdot),
$$
$$ 
\kappa_{ijk}(u, \cdot) : = -\hatl_i(u,\cdot) \cdot D_u\bigl( B \, \hatr_k \bigr) \big (u, \cdot\big ) \, \hatr_j(u, \cdot). 
$$  
The system under study has the form 
$$ 
a_i' - {\mu_i(u, \cdot) \over \eps} \, a_i + O(\eta) \, \sum_j |a_j|
= O(1) \, \sum_{j,k} |a_j| \, |a_k|. 
$$

In a central part of our argument we study the homogeneous system 
\be
\varphi_i ' - {\mu_i(u, \cdot) \over \eps} \, \varphi_i 
+ 
\sum_j \pi_{ij}(u, \cdot)\, \varphi_j = 0, 
\qquad 
\varphi = \bigl( \varphi_1, \ldots, \varphi_N \bigr), 
\label{homogeneous}
\ee
and establish that it has solutions $\varphi_j$, referred to as the {\sl linearized $i$-wave measures} associated with the function $u$, 
which are ``close'' (in a sense to be specified) 
to the  following normalized solutions of the corresponding uncoupled system (obtain by setting $\eta=0$) 
$$ 
\varphi_i^\star := {e^{ - g_i/\eps }  \over I_i}, 
\qquad 
I_i := \int_{-L}^{L}e^{ - g_i / \eps} \, dy, 
\quad 
g_i(y) := - \int_{\rho_i}^y \mu_i(u(x), x) \, dx.  
$$
Here, the constants $\rho_i$ are determined so that the functions $g_i$ are non-negative.  

\

\begin{theorem}  
The system (\ref{homogeneous}) admits a solution 
$\varphi$ such that 
for all $i=1, \ldots, N$ and $y \in [-L,L]$ 
$$ 
\bigl(1 - O(\eta) \big) \, \varphi_i^\star(y)  -  \eps \, O(\eta) \, \sum_j \varphi_j^\star(y) 
\leq \varphi_i(y) 
\leq \bigl( 1 + O(\eta) \bigr) \, \varphi_i^\star(y) 
+ \eps \, O(\eta) \, \sum_j \varphi_j^\star(y). 
$$  
\end{theorem} 

\ 

In constrast with the functions $\varphi_i^\star$, the functions $\varphi_i$ need not be positive.  
Next, to control the total variation of the solutions of (\ref{1.3}), 
we derive Glimm-like estimates on the {\sl wave interaction coefficients} 
$$
F_{ijk}^\star(y) : = \varphi_i^\star(y) \int_{c_i}^y 
        {\varphi_j^\star  \, \varphi_k^\star \over \varphi_i^\star} \, dx 
$$
for some constants $c_i \in [\bLam_i, \Lamb_i]$. We gain useful information on the possible growth of the 
total variation of solutions. Roughly speaking, the coefficient 
$F_{ijk}$ bounds the contribution to the $i$-th family due to 
interactions between waves of the $j$-th and $k$-th characteristic families.


\section{Generalizations}
\label{3-0}

The results above have been also extended 
to relaxation approximations and boundary-value problems. 
In particular, we can handle relaxation approximations associated with
the conservative system (\ref{1.5}) 
\be
\del_t u^\eps + \del_x v^\eps = 0,
\qquad  
\del_t v^\eps + a^2 \, B(u) \del_x u^\eps 
      = {1 \over \eps \, t} \, \bigl( f(u^\eps) - v^\eps \bigr), 
\label{1.9}
\ee 
where $u^\eps = u^\eps(x,t)$ and $v^\eps = v^\eps(x,t)$ are the unknowns, and 
$\eps >0$ is a relaxation parameter.

We also study (\ref{1.3}) and (\ref{1.9}) in the presence of a boundary, when there exists an index $p$ such that
$0 < \Lamb_p$, and that at most one wave family is characteristic, that is, $0 \in (\bLam_p, \Lamb_p)$. 
We consider (\ref{1.1}) on the interval $y \in [0, L]$, and prove the existence of a solution with 
uniformly bounded variation. To handle the boundary layer, we modify the previous definition of the functions
$\varphi_j^\star, j\leq p$ and carefully estimate 
the coefficients $F_{ijk}^\star(y)$ when $\eps \rightarrow 0$.

In addition, following pioneering work by Fan and Slemrod \cite{FanSlemrod}  
who studied the effect of artificial viscosity terms, 
we consider a system arising in liquid-vapor phase dynamics
with {\sl physical} viscosity and capillarity effects taken into account. 
We establish uniform total variation bounds, allowing us to deduce new existence results.  
Our analysis cover both the hyperbolic and the hyperbolic-elliptic regimes and apply to arbitrarily large 
Riemann data. The proofs rely on a new technique of reduction to two coupled scalar equations
associated with the two wave fans of the system.  Strong $L^1$ convergence to a weak solution of bounded variation 
is established in the hyperbolic regime, while in the hyperbolic-elliptic regime
a stationary singularity near the axis separating the two wave fans, or more generally  
an almost-stationary oscillating wave pattern (of thickness depending upon the capillarity-viscosity 
ratio) are observed which prevent the solution to have globally bounded variation.  

\section*{Acknowledgements}
PLF was partially supported by the A.N.R. (Agence Nationale de la Recherche) grant 06-2-134423 (MATH-GR). 
KTJ and PGL were partially supported by a grant (Number 2601-2)
 from the Indo-French Centre 
for the Promotion of Advanced Research, IFCPAR (CEFIPRA).

\end{document}